\documentclass{amsart}

\newcommand{\Z}{\ensuremath{\mathbf Z}}
\newcommand{\Cay}{\ensuremath{\text{Cayley}}}
\newcommand{\outdeg}{\ensuremath{\text{outdeg}}}
\newcommand{\indeg}{\ensuremath{\text{indeg}}}
\newcommand{\Aut}{\ensuremath{\text{Aut}}}

\newtheorem{theorem}{Theorem}
\newtheorem{lemma}{Lemma}
\newtheorem{corollary}{Corollary}
\newcommand{\bt}{\begin{theorem}}
\newcommand{\et}{\end{theorem}}
\newcommand{\bl}{\begin{lemma}}
\newcommand{\el}{\end{lemma}}
\newcommand{\bc}{\begin{corollary}}
\newcommand{\ec}{\end{corollary}}
\newcommand{\beq}{\begin{equation}}
\newcommand{\eeq}{\end{equation}}
\newcommand{\benum}{\begin{enumerate}}
\newcommand{\eenum}{\end{enumerate}}

\begin{document}
\title[The Caccetta-H{\"a}ggkvist conjecture and additive number theory]{The Caccetta-H{\"a}ggkvist conjecture\\ and additive number theory}
\author{Melvyn B. Nathanson}
\address{Department of Mathematics\\ 
Lehman College (CUNY)\\ 
Bronx, New York 10468}
\email{melvyn.nathanson@lehman.cuny.edu}
\thanks{This is the text of two talks in the New York Number Theory Seminar on February 9 and February 16, 2006.  These notes are based on lectures by Matt DeVos and others at the American Institute of Mathematics (AIM) workshop on the Caccetta-H{\"a}ggkvist conjecture in Palo Alto on Janaury 30-February 3, 2006, and on discussions with several conference participants.  I wish to thank AIM for sponsoring this workshop.}
\thanks{The work of M.B.N. is supported in part by grants from the NSA Mathematical Sciences Program and the PSC-CUNY Research Award Program.}

\begin{abstract}
The Caccetta-H{\"a}ggkvist conjecture states that if $G$ is a finite directed graph with at least $n/k$ edges going out of each vertex, then $G$ contains a directed cycle of length at most $k$.  Hamidoune used methods and results from additive number theory to prove the conjecture for Cayley graphs and for vertex-transitive graphs.  This expository paper contains a survey of results on the Caccetta-H{\"a}ggkvist conjecture, and complete proofs of the conjecture in the case of Cayley and vertex-transitive graphs.
\end{abstract}

\maketitle

\section{Many edges imply short cycles}
A \emph{finite directed graph} \index{directed graph} \index{graph!directed} $G = (V,E)$ consists of a finite set $V=V(G)$ of vertices and a finite set $E=E(G)$ of edges, where an edge $e = (v,v')$ is an ordered pair of vertices.  If $e = (v,v')$ is an edge, then the vertex $v$ is called the \emph{tail} of $e$, and $v'$ is called the \emph{head} of $e$.
The {\em outdegree} \index{outdegree} of a vertex $v \in V$, denoted $\outdeg_G(v)$, is the number of edges $e \in E$ of the form $(v,v')$, that is the number of edges with tail $v$.
The {\em indegree} \index{indegree} of a vertex $v' \in V$, denoted $\indeg_G(v')$,  is the number of edges $e \in E$ of the form $(v,v')$, that is the number of edges with head $v'$.

Let $v$ and $v'$ be distinct vertices of the finite directed graph $G$.
A \emph{directed path of length ${\ell}$} in $G$ from vertex $v$ to vertex $v'$ is a sequence of ${\ell}$ edges 
 \[
 (v_0,v_1), (v_1,v_2),\ldots, (v_{{\ell}-1},v_{\ell})
 \]
 such that $v=v_0$ and  $v' = v_{\ell}$. A \emph{directed cycle of length ${\ell}$} in $G$ is a sequence of ${\ell}$ edges 
 $
 (v_0,v_1), (v_1,v_2),\ldots, (v_{{\ell}-1},v_{\ell})
 $
 such that $v_0 = v_{\ell}$.  
 A \emph{loop} is a cycle of length 1, that is, an edge of the form $(v,v)$.  A cycle of length 2 is called a \emph{digon}, and consists of two edges of the form $(v_0,v_1)$ and $(v_1,v_0),$  where $v_0\neq v_1.$  A {\em directed triangle} is a cycle of length 3 of the form $(v_0,v_1), (v_1,v_2), (v_2,v_0)$, where the vertices $v_0,v_1,v_2$ are distinct.

It is reasonable to expect that a finite directed graph with many edges should have many cycles, and, in particular, should have short cycles. 
A quantitative expression of this intuition is the following:
\begin{quote}
There is a function $f(r)$ such that if $G$ is a finite directed graph with $n$ vertices and if there are at least $r$ edges going out of every vertex in $G$, then $G$ contains a cycle of length at most $f(r)n$.  
\end{quote}

A theorem of Chv{\' a}tal-Szemer{\' e}di~\cite{chva-szem83} shows that this is true with $f(r) = 2/(r+1)$.
We start with a simple averaging argument.

\bl   \label{CH:lemma:average}
Let $G = (V,E)$ be a finite directed graph such that $\outdeg_G(v) \geq r$ for every $v \in V$.  There exists a vertex $v_0 \in V$ such that $indeg_G(v) \geq r$.
\el

 \begin{proof}
Suppose that every vertex of $G$ has outdegree at least $r$.  The number of edges $|E|$ satisfies the inequality
\[
 |V|r \leq \sum_{v\in V} \outdeg_G(v)= |E| = \sum_{v\in V} \indeg_G(v) \leq |V|\max\{\indeg_G(v) : v\in V \}
\]
and so there exists a vertex $v_0 \in V$ such that $\indeg_G(v_0) \geq r$.
 \end{proof}

\bt[Chv{\' a}tal-Szemer{\' e}di]   \label{CH:theorem:CSupperbound}
Let $r$ be a positive integer.  If $G = (V,E)$ is a finite directed graph with $|V| = n \geq r$ vertices such that $\outdeg_G(v) \geq r$ for all $v \in V$, then $G$ contains a cycle of length at most $2n/(r+1)$
\et

\begin{proof}
The proof is by induction on $n$.  
If $n=r$, then $G$ is the complete directed graph on $n$ vertices and contains a loop at each vertex.  Let $n \geq r+1$ and assume that the Theorem holds for all graphs with less than $n$ vertices.
The number $|E|$ of edges in the graph satisfies the inequality
\[
r|V| \leq  \sum_{v\in V} \outdeg_G(v) = |E| = \sum_{v\in V} \indeg_G(v)
\]
and so there is a vertex $v_0 \in V$ such that $\indeg_G(v_0) \geq r$.
Let $A$ denote the set of vertices $a \in V$ such that $(a,v_0) \in E$ and 
let $B$ denote the set of vertices $b \in V$ such that $(v_0,b) \in E$.
If $v_0 \in A\cup B$, then $(v_0,v_0)$ is an edge and so $G$ contains a loop, that is, a cycle of length 1.   Similarly, if $v_0 \notin A\cup B$ and $A \cap B \neq \emptyset$, then there is a vertex $v \in V$ such that $(v_0,v)$ and $(v,v_0)$ are both edges in $G$ and so $G$ contains a digon, that is, a cycle of length $2 \leq 2n/(r+1)$.  Therefore, we can assume that the sets $A$, $B$, and $\{v_0\}$ are pairwise disjoint.  

Since $|A|\geq r$ and $|B| \geq r$, it follows that $n \geq 2r+1$ and 
\[
\frac{2n}{r+1} \geq \frac{4r+2}{r+1} \geq 3.
\]
Let $b \in B$ and let $(b,v)$ be an edge in $E$.  If $v =a \in A$, then $(v_0,b)$, $(b,a)$, and $(a,v_0)$ is a directed triangle in  $G$, that is, a cycle of length $3 \leq 2n/(r+1)$.  Therefore, we can also assume that $v\notin A$ for every edge $(b,v) \in E$.

Let $v \in V \setminus (A \cup \{v_0\})$.  Let  $d_A$ denote the number of edges in $E$ of the form $(v,a)$ with $a \in A$, let $d_B$ denote the number of edges in $E$ of the form $(v,b)$ with $b \in B$, and let $d_C$ denote the number of edges in $E$ of the form $(v,v')$ with $v' \notin A\cup B$.  Since $(v,v_0) \notin E$, it follows that
\[
d_A + d_B + d_C = \outdeg_G(v) \geq r.
\]
Let $d' = \min(d_A,|B|-d_B)$.  Choose vertices $b_1,b_2,\ldots,b_{d'} \in B$ such that $(v,b_i) \notin E$ for $i =1,\ldots, d'$, and let 
\[
E'(v) = \{ (v,b_i) : i =1,\ldots, d'\}.
\]
An ordered pair in the set $E'(v)$ will be called a "new edge."
Note that if $v\in B$, then $d_A = 0$ and $E'(v) = \emptyset$.

We construct a new graph $G' = (V',E')$ as follows.  Let 
\[
V' = V \setminus (A \cup \{v_0\}.
\]
Then
\[
n' = |V'| = |V| - |A| - 1 \leq n-r-1.
\]
Let
\[
E'_0 = \{ (v,v') \in E : v,v' \in V'\}
\]
and
\[
E' = E'_0 \cup \bigcup_{v\in V'} E'(v).
\]
If $b \in B$, then $\outdeg_{G'}(b) = \outdeg_G(v) \geq r$.
If $v \in V' \setminus B$, then $\outdeg_{G'}(v) = d' + d_B + d_C$.  
If $d' = d_A$, then $\outdeg_{G'}(v) = d_A + d_B + d_C  \geq r$.  
If $d' = |B|-d_B$, then $\outdeg_{G'}(v) = |B| + d_C \geq |B| \geq r$.  
Therefore, every vertex in $G'$ has outdegree at least $r$.  
Since $|V'| \leq n-r-1,$ the induction hypothesis implies that $G'$ contains a cycle $\mathcal{C}'$ of length 
\[
\ell' \leq \frac{2|V'|}{r+1} \leq \frac{2(n-r-1)}{r+1}. 
\]
If $(v,b)$ is a "new edge" in this cycle, that is, if $(v,b) \in E(v)$, then there exists $a \in A$ such that $(v,a)$ is an edge in $E$, and
\beq  \label{CH:newedge}
(v,a), (a,v_0),(v_0,b)
\eeq
is a directed path in $E$.    Suppose that $\mathcal{C}'$ contains exactly $m$ new edges.  Replacing every new edge $(v,b)$ in the cycle $\mathcal{C}'$ with three old edges of the form~\eqref{CH:newedge}, we obtain a cycle $\mathcal{C}$ in the original graph $G$ of length 
\[
\ell' + 2m \leq \frac{2(n-r-1)}{r+1} + 2m.
\]
The vertex $v_0$ occurs exactly $m$ times in this cycle, and so the cycle decomposes into $m$ cycles, and the sum of the lengths of these $m$ cycles is exactly $\ell' + 2m$.  This implies that $G$ contains a cycle of length at most 
\begin{align*}
\frac{\ell' + 2m}{m} 
& \leq \left(\frac{1}{m} \right)\left( \frac{2(n-r-1)}{r+1} + 2m\right)  \\
& =  \frac{2n}{m(r+1)} - \frac{2}{m} + 2 \\
& \leq  \frac{2n}{r+1}.
\end{align*}
This completes the proof.
\end{proof}

For every real number $t$, let $\lceil t \rceil$ denote the smallest integer $n \geq t.$ 

Shen~\cite{shen00} obtained a significant improvement of Theorem~\ref{CH:theorem:CSupperbound}.
He proved that if $G = (V,E)$ is a finite directed graph with $|V| = n \geq r$ vertices such that $\outdeg_G(v) \geq r$ for all $v \in V$, then $G$ contains a cycle of length at most
\[
3\left\lceil \left( \ln \frac{2+\sqrt{7}}{3}  \right) \frac{n}{r} \right\rceil \approx \frac{1.312n}{r}.
\]

Caccetta and H{\"a}ggkvist~\cite{cacc-hagg78} \index{Caccetta-H{\"a}ggkvist conjecture} made a strong assertion about the existence of short cycles in directed graphs with many edges.  Their conjecture states:
\begin{quote}
If $G$ is a finite directed graph with $n$ vertices such that every vertex has outdegree at least $r$, then the graph contains a directed cycle of length at most $\lceil n/r \rceil$.
\end{quote}
The \emph{girth} \index{girth} of a graph is the length of the shortest cycle in the graph.  We can restate the Caccetta-H{\"a}ggkvist conjecture as follows:
If every vertex in a finite directed graph has outdegree at least $r$, then the girth of the graph is at most $\lceil n/r \rceil$. 

If $G$ is a finite directed graph such that every vertex is the tail of at least one edge, that is, if $\outdeg_G(v) \geq 1$ for all $v \in V$, then $G$ contains a directed cycle.  This is the case $r=1$ of the Caccetta-H{\"a}ggkvist conjecture.  The conjecture has been proved for $r =2$ by Caccetta, and H{\"a}ggkvist~\cite{cacc-hagg78}, for $r = 3$ by Hamidoune~\cite{hami87a}, and for $r = 4$ and $5$ by Ho{\' a}ng and Reed~\cite{hoan-reed87}.  Shen~\cite{shen00} proved that conjecture holds for all $r \geq 2$ and $n \geq 2r^2-3r+1$.   The conjecture has also been proved ``up to an additive constant'' in the following form:  If $G$ is a finite directed graph with $n$ vertices such that every vertex has outdegree at least $r$, then  the girth of $G$ is at most $\lceil n/r \rceil + c$.   Chv{\' a}tal and Szemer{\' e}di~\cite{chva-szem83} obtained $c = 2500$, Nishimura~\cite{nish88} obtained $c = 304$, and 
Shen~\cite{shen00} obtained $c = 73$.

 The following example shows that the upper bound in the Caccetta-H{\"a}ggkvist conjecture is best possible.

\bt[Behzad, Chartrand, and Wall~\cite{behz-char-wall70}]   
\label{CH:theorem:lowerbound}
Let $r$ be a positive integer.  For every integer $n \geq r$ there is a graph $G = (V,E)$ with $|V| = n$ vertices such that $\outdeg_G(v) \geq r$ for all $v \in V$ and the girth of $G$ is exactly $\lceil n/r \rceil$.
\et

\begin{proof}
Let $n \geq r$ and $A = \{1,2,\ldots,r\}$.  Consider  the additive group $\Z/n\Z = \{ \overline{0}, \overline{1}, \overline{2}, \ldots, \overline{n-1}\}$, where $\overline{x} = x + n\Z$.  Let $G = (V,E)$ be the graph whose vertices are the congruence classes in $\Z/n\Z$  and whose edges are the ordered pairs of the form $(\overline{x}, \overline{x}+\overline{a})$, where $\overline{x} \in \Z/n\Z$ and $a \in A$.  Let
\beq   \label{CH:cycle-ell}
(v_0,v_1), (v_1,v_2),\ldots, (v_{{\ell}-1},v_{\ell})
\eeq
be a cycle of length $\ell$, where $v_0 = v_{\ell}$ and 
$v_i = v_{i-1} + \overline{a_i}$ for $a_i \in A$ and $i = 1,\ldots, \ell.$
Then $\overline{a_1} + \overline{a_2} + \cdots + \overline{a_{\ell}} = 0$ in the group $\Z/n\Z$, and so 
\[
a_1 + a_2 + \cdots + a_{\ell} \equiv 0 \pmod{n}.
\]
Since
\[
0 <  \ell \leq a_1 + a_2 + \cdots + a_{\ell} \leq r\ell
\]
it follows that $r\ell \geq n$ and so $\ell \geq n/r$.  Therefore, $\ell \geq \lceil n/r \rceil$ and the girth of the graph $G$ is at least $\lceil n/r \rceil$.

Conversely, let $n = \ell r-s$, where $\ell$ and $s$ are integers and $0 \leq s \leq r-1$.  Then $\ell = \lceil n/r \rceil \geq 1$.  If $s=0$, let $a_i = r$ for $i=1,\ldots, \ell$.  If $1 \leq s\leq r-1$, let  $a_i = r$ for $i=1,\ldots, \ell -1$ and $a_{\ell} = r-s$.  Then $a_i \in A$ for $i=1,\ldots, \ell$ and $a_1 + a_2 + \cdots + a_{\ell} = n$.  For any $\overline{x} \in V = \Z/n\Z$, let
\[
v_i = \overline{x} + \overline{a_1}  + \overline{a_2}  \cdots  + \overline{a_i}
\]
for $i = 0,1,\ldots, \ell$.  Then~\eqref{CH:cycle-ell} is a cycle in $G$, and the girth of $G$ is at most $\lceil n/r\rceil$.  This completes the proof.
\end{proof}

\subsection*{Exercises}
\benum
\item
Let $G = (V,E)$ be a directed graph.  Prove that if $\outdeg_G(v) = d^+$ for every $v\in V$ and  if $\indeg_G(v) = d^-$ for every $v\in V$,  then $d^+=d^-$.
The graph $G$ is called \emph{regular of degree $d$} if $\outdeg_G(v)=\indeg_G(v) = d$ for every $v \in V$.\index{regular graph}\index{graph!regular}

\item
The directed graph $G=(V,E)$ is \emph{path-connected} \index{path-connected} if for every pair of distinct vertices $v,v' \in V$ there is a directed path from $v$ to $v'$.   If $G$ is path-connected, then the \emph{distance}\index{graph!distance} from vertex $v$ to vertex $v'$ is the length of the shortest directed path from $v$ to $v'$.  The \emph{diameter}\index{diameter}\index{graph!diameter} of a path-connected graph $G$ is the maximum distance between two vertices of $G$.  Prove that if $G$ is a path-connected directed graph with diameter $D$ and girth $g$, then $g \leq D+1$.

\item
Let $G = (V,E)$ be a finite directed graph with neither loops nor digons.  For every vertex $v \in V$, the \emph{first neighborhood} \index{first neighborhood}
$N^+(v)$ \index{neighborhood!first}  consists of all vertices $v'$ such that $(v,v') \in E$.  The \emph{second neghborhood}\index{second neighborhood}\index{neighborhood!second} $N^{++}(v)$ is the set of all vertices $v'' \in V$ such that (i) there is a vertex $v'\in N^+(v)$ with $(v,v') \in E$ and $(v',v'') \in E$, and (ii) $v'' \notin N^+(v)$.  Seymour's \emph{second neighborhood conjecture} \index{second neighborhood conjecture} states that there is a vertex $v\in V$ such that 
\[
|N^+(v)| \leq |N^{++}(v)|.
\]
Show that the second neighborhood conjecture implies that  if $G$ is a graph with $n$ vertices such that (i) $G$ contains no loops and no digons and (ii) every vertex of $G$ has indegree and outdegree at least $n/3$, then the Caccetta-H\" aggkvist conjecture is true for $G$, that is, $G$ contains a directed triangle.
(Hint: Prove that there is a vertex $v \in V$ such that $(v'',v)\in E$ for some $v'' \in N^{++}(v)$.)

\eenum

\section{Directed triangles in directed graphs}

An equivalent form of the Caccetta-H{\"a}ggkvist conjecture is the following:  If $G$ is a finite directed graph with $n$ vertices such that every vertex has outdegree at least $n/k$, then $G$ contains a directed cycle of length at most $k$.  If $k=1$, then every vertex has degree $n$, so the graph contains loops, which are cycles of length 1.  

\bt[Caccetta, and H{\"a}ggkvist~\cite{cacc-hagg78}]
If $G$ is a finite directed graph with $n$ vertices such that every vertex has outdegree at least $n/2$, then $G$ contains a loop or a digon, that is, a directed cycle of length at most 2.
\et

\begin{proof}
Suppose that every vertex of $G$ has outdegree at least $n/2$.  By Lemma~\ref{CH:lemma:average}, there exists a vertex $v_0 \in V$ such that $\indeg_G(v_0) \geq n/2$.
Let $V' = \{v'\in V : (v',v_0) \in E\}$ and let $V'' = \{v'' \in V : (v_0,v'') \in E\}$.
Since $|V'| \geq n/2$ and $|V''| \geq n/2$, it follows that the sets $V'$, $V''$, and $\{v_0\}$ cannot be pairwise disjoint.  If $v_0 \in V' \cup V''$, then $G$ contains a loop.  Otherwise, $V' \cap V'' \neq \emptyset$, and $G$ contains a digon.  
\end{proof}

For $k = 3$, the Caccetta-H{\"a}ggkvist conjecture asserts that if $G$ has outdegree at least $n/3,$ then $G$ contains a cycle of length at most 3, that is, a loop, digon, or triangle.  This is a famous unsolved problem in graph theory.  

\bt[Caccetta, and H{\"a}ggkvist~\cite{cacc-hagg78}]
 \label{CH:theorem:triangle}
Let 
\[
c_0 = \frac{3-\sqrt{5}}{2}\approx 0.3820 \ldots.
\]
If $G = (V,E)$ is a finite directed graph with $|V| = n$ vertices such that $\outdeg_G(v) \geq c_0 n$ for all $v \in V$, then $G$ contains a cycle of length at most 3.
\et

\begin{proof}
Let $0 < c < 1$, and let $G = (V,E)$ be a directed graph with $n$ vertices such that $\outdeg_G(v) \geq cn$ for all $v \in V$ and $G$ does not contain a loop, digon, or triangle.  We shall prove that $c < (3-\sqrt{5})/2$.

By Lemma~\ref{CH:lemma:average}, the graph $G$ contains a vertex $v_0$ such that $\indeg_G(v_0) \geq cn$.  Let $A$ be the set of vertices $a$ such that $(a,v_0) \in E$, and let $B$ be the set of vertices $b$ such that $(v_0,b) \in E$.  Then $|A| \geq cn$ and $|B| \geq cn$.   If $v_0 \in A\cup B$, then $G$ contains a loop.  Similarly, if $v_0 \notin A\cup B$ and $A \cap B \neq \emptyset$, then $G$ contains a digon.  Therefore, we can assume that the sets $A$, $B$ and $\{v_0\}$ are pairwise disjoint.

Let $G'$ be the complete subgraph of $G$ induced by $B$, that is, $G'$ is the graph whose vertex set is $B$ and whose edges are all ordered pairs $(b,b')$ such that $b,b' \in B$ and $(b,b') \in E$.  Since $|B|<n$, it follows from the induction hypothesis that if $\outdeg_{G'}(b) \geq c|B|$ for all $b\in B$, then the graph $G'$ contains a triangle, and so $G$ contains a triangle.  Therefore, we can assume that there is a vertex $b_0 \in B$ such that $\outdeg_{G'}(b_0) < c|B|$.  Let $W$ be the set of all vertices $w \in V\setminus B$ such that $(b_0,w) \in E$.  Since $\outdeg_G(b_0) \geq cn$, it follows that 
\[
|W| = \outdeg_G(b_0) - \outdeg_{G'}(b_0) >  cn -  c|B|.
\] 
If $v_0 \in W$, then $G$ contains a digon.  If $A \cap W \neq \emptyset$, then $G$ contains a triangle.  Therefore, we can assume that the sets $A$, $B$, $W$, and $\{v_0\}$ are pairwise disjoint subsets of $V$, and so
\[
n \geq |A|+ |B| + |W| + 1 > 2cn + (1-c)|B| + 1 > 3cn - c^2.
\]
This implies that
\[
c^2-3c+1 > 0
\]
and so 
\[
c < \frac{3-\sqrt{5}}{2}.
\]
Therefore, if $c \geq (3-\sqrt{5})/2$, then $G$ contains a cycle of length at most 3.
\end{proof}

The constant $c$ in Theorem~\ref{CH:theorem:triangle} has been reduced by Bondy~\cite{bond97}, who obtained
\[
c_0 = \frac{2\sqrt{6}-3}{5} \approx 0.3798
\]
and by Shen~\cite{shen98}, who obtained 
\[
c_0 = 3-\sqrt{7} = 0.3542\ldots.
\]

\section{Kemperman's theorem for nonabelian groups}
In the following sections we shall prove the Caccetta-H{\"a}ggkvist conjecture for two important classes of finite directed graphs:  Cayley graphs and vertex-transitive graphs.  The proof depends on a result of Kemperman that gives an upper bound for the growth of certain subsets of a group.

Let $\Gamma$ be a finite group, written multiplicatively, and let $(A,B)$ be a pair of finite subsets of $\Gamma$.  The \emph{product set} $A\cdot B$ is the set
\[
A\cdot B = \{ab : a\in A \text{ and } b\in B\}.
\]
We define the iterated product sets $B^2 = B\cdot B$ and $B^k = B\cdot B^{k-1}$ for all $k \geq 2.$  Then
\[
B^k = \{b_1b_2\cdots b_k : b_i \in B \text{ for $i = 1,2,\ldots,k$}\}.
\]
We usually write $AB$ instead of $A\cdot B$.
For $x \in \Gamma$ and $A \subseteq \Gamma$, let
\[
Ax = A \{x\} = \{ ax: a \in A\}
\] 
and
\[
xA =  \{x\}A = \{xa : a \in A\}.
\]
Let $|X|$ denote the cardinality of the set $X$.  For any pair $(A,B)$ of finite subsets of $\Gamma$, we define
\[
k(A,B) = |A|+|B|-|A+B|.
\]

\begin{theorem}[Kemperman~\cite{kemp56a}]  \label{CH:theorem:kemperman}
Let $\Gamma$ be a finite group and let $(A,B)$ be a pair of finite subsets of $\Gamma$ such that 
\begin{enumerate}
\item[(i)]
$1 \in A \cap B$ 
\item[(ii)]
If $a\in A$, $b\in B$, and $ab=1$, then $a=b=1.$  
\end{enumerate}
Then
\[
|AB| \geq |A| + |B| -1. 
\]
\end{theorem}

\begin{proof}
Suppose there exist pairs $(A,B)$ of subsets of $\Gamma$ such that $A$ and $B$ satisfy conditions~(i) and~(ii), but $|AB| < |A| + |B| - 1$.  Equivalently, 
\begin{equation}  \label{CH:kempineq}
k(A,B) = |A|+|B|-|AB| \geq 2.
\end{equation}
Consider pairs $(A,B)$ that have the maximum value of $k(A,B)$.  Among all such pairs, choose $(A,B)$ with the minimum value of $|A|$.

Since $1 \in A \cap B$, it follows that $AB \supseteq A \cup B$.  If $A \cap B = \{1\}$, then $|A B| \geq |A \cup B| = |A| + |B|-1$, which contradicts~\eqref{CH:kempineq}.  Therefore, there exists $x \in A \cap B$ with $x \neq 1$.
We introduce the sets
\begin{align*}
A_1 & = Ax^{-1} \cap A \\
B_1 & = xB \cup B \\
A_2 & = Ax \cup A \\
B_2 & = x^{-1}B \cap B .
\end{align*}
Then
\[
A_1B_1 \subseteq AB \text{ and } A_2B_2 \subseteq AB .
\]

We shall show that $A_1$ is a proper subset of $A$.  Note that $x \in A$.  If $x^k \in A$ for all positive integers $k$, then $x$ has finite order $m$ since $A$ is finite, and so $x^{m-1} = x^{-1} \in A$.  Since $x \in B$, we have $1 = x^{-1} \cdot x \in AB$  and so $x=1$, which is a contradiction.  It follows that there must exist a largest positive integer $k$ such that $x^k \in A$.  If $x^k \in A_1$, then $x^k \in Ax^{-1}$ and so $x^{k+1} \in A$, which is a contradiction.  Therefore,
$x^k \in A \backslash A_1$ and $|A_1| < |A|$.

By Exercise~\ref{CH:exer:Ax},
\[
A_1x = (Ax^{-1} \cap A)x = A \cap Ax
\]
and so
\[
|A_1| = |A_1x| = |A\cap Ax|
\]
and
\beq   \label{CH:4}
|A_1|+|A_2| = |A\cap Ax| + |Ax \cup A| = |A|+|Ax| = 2|A|.
\eeq
Similarly,
\beq   \label{CH:5}
|B_1|+|B_2| = 2|B|.
\eeq
Adding~\eqref{CH:4} and~\eqref{CH:5}, we obtain
\beq  \label{CH:6}
(|A_1|+|B_1|)+(|A_2| + |B_2|)= 2( |A| + |B|).
\eeq

We shall show that the pairs $(A_1 , B_1 )$ and $(A_2 , B_2 )$ also satisfy conditions~(i) and~(ii).
Since $1 \in A\cap B$ and $x \in A \cap B$, it follows that
\[
1 \in A_1\cap B_1 \text{ and } 1 \in A_2 \cap B_2
\]
and so the first condition is satisfied.

Suppose that  $a_1\in A_1$, $b_1 \in B_1$, and $a_1b_1 = 1$.  Then $b_1 \in B$ or $b_1 \in xB$.  
If $b_1 \in B$, then $a_1\in A_1 = Ax^{-1} \cap A \subseteq A$ implies that $a_1=b_1=1$.  
On the other hand, if $b_1 \in xB$, then $b_1 = xb$ for some $b \in B$.  Since $a_1 \in A_1 \subseteq Ax^{-1}$, there exists $a\in A$ such that $a_1=ax^{-1}$.  Then $1 = a_1b_1 = ax^{-1}xb = ab$, and so $a=b=1$.  It follows that $a_1 = x^{-1} \in A$, which is impossible because $x\in B$ and $x \neq 1$.   Therefore, the pair $(A_1,B_1)$ satisfies condition~(ii).  Similarly, $(A_2,B_2)$ satisfies condition~(ii).

By the maximality of $k(A,B)$, 
\beq   \label{CH:2}
|A_1|+|B_1|-|A_1B_1| = k(A_1,B_1) \leq k(A,B) = |A|+|B|-|AB|
\eeq
and
\beq   \label{CH:3}
|A_2|+|B_2|-|A_2B_2| = k(A_2,B_2) \leq k(A,B) =  |A|+|B|-|AB|.
\eeq
Adding and rearranging~\eqref{CH:2} and~\eqref{CH:3}, we obtain
\[
2|AB| \leq |A_1B_1| + |A_2B_2|.
\]
Since $A_1B_1 \subseteq AB$, $A_2B_2 \subseteq AB$, we also have
\[
|A_1B_1| \leq |AB| \text{ and } |A_2B_2| \leq |AB|
\]
and so
\[
|A_1B_1| = |A_2B_2| =|AB|.
\]
Therefore,
\[
k(A_1,B_1) + k(A_2,B_2) = 2 k(A,B).
\]
The maximality of $k(A,B)$ implies that
\[
k(A_1,B_1) =  k(A_2,B_2) = k(A,B)
\]
but this is impossible since the inequality $|A_1|<|A|$ contradicts the minimality of $|A|$.
This completes the proof.
\end{proof}

\begin{theorem}  \label{CH:theorem:kemperman-k}
Let $\Gamma$ be a group and let $B$ be a finite subset of $\Gamma$ with $1 \in B$.   If the only solution of the equation $b_1b_2\cdots b_k = 1$ with $b_i \in B$ for all $i=1,\ldots, k$ is $b_1 = b_2 = \cdots = b_k = 1,$ then
\[
|B^k| \geq k|B|-k+1.
\]
\end{theorem}

\begin{proof}
Exercise~\ref{CH:problem:kemperman}.
\end{proof}

\subsection*{Exercises}
\benum
\item \label{CH:exer:Ax}
Let $A$ be a finite subset of a group $\Gamma$ and let $x \in \Gamma$.  Prove that $(Ax^{-1}\cap A)x = A \cap Ax$.

\item  \label{CH:problem:kemperman}
Prove Theorem~\ref{CH:theorem:kemperman-k} by induction on $k$.

\eenum

\section{The Caccetta-H{\"a}ggkvist conjecture for Cayley graphs}

Let $\Gamma$ be a finite group, not necessarily abelian.  We write the group operation multiplicatively.  Let $A$ be a subset of $\Gamma$.  The \emph{Cayley graph} \index{Cayley graph} \index{graph!Cayley} $\Cay(\Gamma,A)$ is the graph whose vertex set is the group $\Gamma,$ and whose edge set consists of all ordered pairs of form $(v,va),$ where $v \in \Gamma$ and $a \in A$.  By Exercise~\ref{CH:exer:inout}, every vertex in $\Cay(\Gamma,A)$ has outdegree $|A|$ and indegree $|A|$.  Moreover, $\Cay(\Gamma,A)$ contains a loop if and only if $1 \in A$, and $\Cay(\Gamma,A)$ contains a digon if and only if $\{a,a^{-1}\} \subseteq A$ for some $a \in \Gamma, a \neq 1.$

\begin{lemma}   \label{CH:lemma:cycle}
Let $\Gamma$ be a finite group and $A \subseteq \Gamma$.  The graph $\Cay(\Gamma,A)$ contains a directed cycle of length $\ell$ if and only if 
$1 \in A^{\ell}$.
\end{lemma}

\begin{proof}
If $1 \in A^{\ell}$, then there exist $a_1,a_2,\ldots,a_{\ell}$ in $A$ such that $a_1a_2 \cdots a_{\ell} = 1$.
For any $v_0 \in \Gamma$, if we define 
\[
v_i = v_{i-1}a_i = v_0a_1a_2\cdots a_i
\]
for $i=1,\ldots,\ell$, then $v_{\ell} = v_0$ and  
\beq      \label{CH:cycle}
(v_0,v_1), (v_1,v_2),\ldots, (v_{{\ell}-1},v_{\ell})
\eeq
is a directed cycle of length ${\ell}$ in $\Cay(\Gamma,A)$.  

Conversely, if~\eqref{CH:cycle} is a directed cycle in $\Gamma$, then $v_0 \in \Gamma$ and there exist $a_1$,\ldots,$a_{\ell} \in A$ such that $v_i=v_{i-1}a_i$ for $i = 1,\ldots,{\ell}.$  This implies that if $1 \leq i < j \leq {\ell}$, then 
\[
v_j = v_{i-1} a_ia_{i+1}\cdots a_j.
\]
In particular,
\[
v_0 = v_{\ell} = v_0a_1a_2\cdots a_{\ell}
\]
and so
\[
1 = a_1\cdots a_{\ell} \in A^{\ell}.
\]
This completes the proof.
\end{proof}

\begin{theorem}[Hamidoune~\cite{hami81a}]  \label{CH:theorem:Cayley}
Let $\Gamma$ be a finite group of order $n$, and let $A$ be a subset of $\Gamma$ such that $|A| \geq n/k.$
Then the graph $\Cay(\Gamma,A)$ contains a cycle of length at most $k$.
\end{theorem}

\begin{proof}
If every cycle in the graph $\Cay(\Gamma,A)$ has length greater than $k$, then Lemma~\ref{CH:lemma:cycle} implies that 
\[
1 \notin A \cup A^2 \cup \cdots \cup A^k.
\]
Let $B = A \cup \{1\}$.  Then $|B|=|A|+1$.  Since the only solution of $1 = b_1b_2\cdots b_k \in B^k$ is $b_1=b_2=\cdots = b_k=1$, it follows from Theorem~\ref{CH:theorem:kemperman-k} that
\[
n =  |\Gamma| \geq |B|^k \geq k|B|-k+1 = k|A|+1 > k|A|
\]
and so $|A| < n/k$, which is false.  Thus, $\Cay(\Gamma,A)$ must contain a cycle of length at most $k$.
\end{proof}

\subsection*{Exercises}
\benum
\item
Prove that the Cayley graph $\Cay(\Gamma,A)$ is path-connected if and only if the semigroup generated by $A$ is $\Gamma.$

\item  \label{CH:exer:inout}
Prove that if $x$ is a vertex in the graph $\Cay(\Gamma,A)$, then $\indeg_G(x) = \outdeg_G(x) = |A|$.

\item
The sequence of edges~\eqref{CH:cycle} is a \emph{simple directed cycle} in $\Gamma$ if $v_i \neq v_j$ for $0 \leq i < j \leq \ell-1$.
Prove that $\Cay(\Gamma,A)$ contains a simple directed cycle of length $\ell$ if and only if there exists a sequence of elements $a_1,\ldots,a_{\ell}$ in $A$ such that $a_ia_{i+1}\cdots a_j = 1$ for $1 \leq i \leq j \leq \ell$ if and only if $i=0$ and $j=\ell$.

\item
Let $\Gamma = \langle x \rangle$  be the cyclic group of order $d(g-1)+1$, written multiplicatively, and let $A = \{x,x^2,\ldots,x^d\}$.  
Prove that the graph $\Cay(\Gamma,A)$ is regular of degree $d$ with girth $g$.
 
\eenum

\section{Graph automorphisms and vertex-transitive graphs}
Let $G=(V,E)$ and $G'=(V',E')$ be finite directed graphs.  The function $x:V \rightarrow V'$ is a \emph{graph isomorphism} \index{graph isomorphism} if $x$ is a bijection and $(v_1,v_2)\in E$ if and only if $(x(v_1),x(v_2))\in E'$.  A \emph{graph automorphism} \index{graph automorphism} is a graph isomorphism from $G$ to $G$.  The automorphisms of a graph form a group, denoted $\Aut(G)$.  We denote the action of an automorphism  $x:V \rightarrow V$ on the vertex $v$ by $xv$.

Let $\Gamma$ be a group of graph automorphisms of $G$, that is, a subgroup of $\Aut(G)$.  The graph $G$ is called \emph{vertex-transitive with respect to $\Gamma$} \index{vertex-transitive graph} \index{graph!vertex-transitive} if, for every pair of vertices $v,v' \in V$, there is an automorphism $x \in \Gamma$ such that $xv = v'$.  We call $G$ \emph{vertex-transitive} if it is vertex-transitive with respect to some group of automorphisms.  In a vertex-transitive graph, $\outdeg_G(v) = \outdeg_G(v')$ for all vertices $v,v' \in V$ (Exercise~\ref{CH:exer:vertexdegree}).  

For example, every Cayley graph is vertex-transitive.  Let $G = \Cay(\Gamma,A)$, where $\Gamma$ is a finite group and $A \subseteq \Gamma$.  To every element $x \in \Gamma$ there is a bijection $x: \Gamma \rightarrow \Gamma$ defined by $v \mapsto xv$ for all $v\in \Gamma$.  If $(v_1,v_2)$ is an edge in $G$, then $v_2 = v_1a$ for some $a\in A$, and $(xv_1,xv_2) = (xv_1,(xv_1)a)$ is also an edge in $G$.  Thus, the map $x \mapsto xv$ is an automorphism of $G$.
In particular, if $v,v' \in \Gamma$ and $x = v'v^{-1}$, then the map $w \mapsto xw=v'v^{-1}w$ sends $v$ to $v'$, and so $\Aut(G)$ acts transitively on $\Gamma$.

We shall prove that the Caccetta-H{\"a}ggkvist conjecture is true for all vertex-transitive graphs.

\begin{theorem}[Hamidoune~\cite{hami81a}]   \label{CH:theorem:vertex}
Let $G$ be a vertex-transitive finite directed graph with $n$ vertices such that $\outdeg_G(v) = d$ for every vertex $v$ of $G$.  Then $G$ contains a cycle of length at most $\lceil n/d \rceil$.
\end{theorem}

\begin{proof}
Let $\Gamma$ be a group of automorphisms that acts transitively on the set $V$ of vertices of a finite directed graph $G$.  For every vertex $v \in V$, the \emph{stabilizer of $v$} \index{stabilizer} is the set 
\[
H_v = \{x\in \Gamma : xv = v\}.
\]
$H_v$ is a subgroup of $\Gamma$.  

Choose a vertex $v_0 \in V$, and let $H_0 = H_{v_0}$.
Since $G$ is vertex-transitive, there is a set $\{x_v\}_{v\in V}$ contained in $\Gamma$ such that $x_v v_0 = v$ for all $v \in V$.  Then
\begin{align*}
H_v  & =  \{x\in \Gamma : xv = v\} \\
& =  \{x\in \Gamma : xx_v v_0 = x_v v_0\} \\
& =  \{x\in \Gamma : x_v^{-1} x x_v v_0 = v_0\} \\
& =  \{x\in \Gamma : x_v^{-1} x x_v \in H_0\} \\
& = x_v H_0 x_v^{-1}.
\end{align*}
The subgroup $H_0$ is normal in $\Gamma$ if and only if $H_v = H_0$ for all $v \in V$.

For all $x\in \Gamma$, we have $xv_0 = v$ if and only if $xv_0 = x_v v_0$ if and only if  $x_v^{-1}x \in H_0$
if and only if $x \in x_v H_0$.    Therefore,
\[
x_vH_0 = xH_0 = \{ x \in \Gamma : xv_0 = v\}.
\]
Let
\[
\Gamma_0 = \Gamma/H_0 = \{xH_0 : x\in \Gamma\}
\]
denote the set of left cosets of $H_0$.  The map $\phi: V \rightarrow \Gamma/H_0$ defined by $v \mapsto x_v H_0$ is a one-to-one correspondence between the vertices of $G$ and the left cosets of $H_0$, and so
\[
|\Gamma_0| = |\Gamma/H_0| = |V| = n
\]
and
\[
|\Gamma| = |\Gamma/H_0| |H_0| = n|H_0|.
\]

We can use the left cosets of $H_0$ to describe the edges in the graph $G$.  Let
\[
A = \{x \in \Gamma:  (v_0, xv_0) \in E \} .
\]
If $x \in A$ and $h \in H_0$, then 
$(v_0,xhv_0) = (v_0,xv_0) \in E$ and so $xh \in A$, hence $xH_0 \subseteq A$.  It follows that $A$ is a union of left cosets of $H_0$.  Let
\[
A_0 = \{ x H_0 : xH_0 \subseteq A \} = \{ x H_0 : (v_0, xv_0) \in E \} = \{ x_v H_0 : (v_0, v) \in E \}.
\]
Then
\[
|A_0| = \outdeg_G(v_0) = d
\]
and
\[
|A| = |A_0| |H_0| = d|H_0|.
\]
Since $\Gamma$ is a group of automorphisms of the graph $G$, the ordered pair $(v,v')$ is an edge of $G$ if and only if 
$
( x_v^{-1} v,x_v^{-1}v')  = (v_0, x_v^{-1} x_{v'} v_0) \in E.
$
Thus,
\[
(v,v' ) \in E \text{ if and only if } x_v^{-1} x_{v'} H_0 \in A_0.
\]

Suppose that  $H_0$ is a normal subgroup of $\Gamma$.  Then $\Gamma_0$ is a group of order $n$.  The graph $\Cay(\Gamma_0,A_0)$ has $|\Gamma_0| = n$ vertices and $|A_0| = d$ edges.  By Theorem~\ref{CH:theorem:Cayley},  $\Cay(\Gamma_0,A_0)$ contains a cycle of length not exceeding  $\lceil n/|A_0|\rceil = \lceil n/d\rceil$.

We shall show that the graphs $G$ and $\Cay(\Gamma_0,A_0)$ are isomorphic.  Recall the bijection $\phi: V \rightarrow \Gamma_0$ defined by $\phi(v) =  x_v H_0$.  Let $(v,v') \in E$, and define $x =  x_v^{-1} x_{v'}$.  Then  $xH_0 = x_v^{-1} x_{v'} H_0 \in A_0$ and $(x_vH_0)(xH_0)=x_{v'}H_0$, hence $(\phi(v),\phi(v')) = (x_vH_0,x_{v'}H_0)$ is an edge in $\Cay(\Gamma_0,A_0)$.   Conversely, if  $(\phi(v),\phi(v')) = (x_vH_0,x_{v'}H_0)$ is an edge in $\Cay(\Gamma_0,A_0)$, then there is a coset $xH_0\in A_0$ such that $x_vH_0 xH_0 = x_vxH_0 = x_{v'}H_0$.  It follows that $x_v^{-1}x_{v'}H_0 \in A_0$ and so $(v,v')\in E$.  Thus, the map $\phi: V \rightarrow \Gamma_0/H_0$ is a graph isomorphism, and so $G$ contains a cycle of length at most $n/d$.

Next we consider the general case when $H_0$ is not necessarily a normal subgroup of $\Gamma$. 
The Cayley graph $\Cay(\Gamma,A)$ contains $|\Gamma| = n|H_0|$ vertices, and and the outdegree of every vertex is $|A| = d|H_0|$.  By Theorem~\ref{CH:theorem:Cayley}, $\Cay(\Gamma,A)$ has a cycle of length $\ell$, where
\[
\ell \leq 
\left\lceil \frac{|\Gamma_0|}{|A|} \right\rceil =
\left\lceil \frac{n|H_0|}{d|H_0|} \right\rceil = \left\lceil\frac{n}{d} \right\rceil.
\]
By Lemma~\ref{CH:lemma:cycle}, there exist elements $a_1,\ldots,a_{\ell} \in A$ such that 
\[
a_1 a_2 \cdots a_{\ell -1 }a_{\ell}  = 1.
\]
Let $v_0 \in \Gamma$ and consider the sequence of vertices $v_0,v_1,\ldots, v_{\ell}$, where
\[
v_i = a_1 a_2 \cdots a_{i} v_0
\]
for $i = 1,2,\ldots,\ell$.  
Then $ (v_0, a_iv_0) \in E$ since $a_i \in A$ for $i = 1,2,\ldots,\ell$.   Since $\Gamma$ is a group of graph automorphisms, we have
\[
(v_{i-1},v_i) =  ( (a_1 a_2 \cdots a_{i-1})v_0, (a_1 a_2 \cdots a_{i-1})a_iv_0) \in E \text{ for $i = 1,2,\ldots,\ell$}
 \]
and
\[
v_{\ell} = a_1 a_2 \cdots a_{\ell -1 }a_{\ell} v_0 = 1\cdot v_0 = v_0.
\]
It follows that 
\[
(v_0,v_1), (v_1,v_2),\ldots, (v_{\ell -1},v_{\ell}) = (v_{\ell -1},v_0)
\]
is a cycle in $G$ of length $\ell \leq \lceil n/d \rceil$.
This completes the proof.
\end{proof}

\subsection*{Exercise}
\benum
\item  \label{CH:exer:vertexdegree}
If $G = (V,E)$ is a vertex-transitive graph, then there is an integer $d$ such that $G$ is regular of degree $d$, that is, $\outdeg_G(v) = \indeg_G(v) = d$ for all $v \in V$.
\eenum

\section{Additive compression}
Let $V_1$ and $V_2$ be finite disjoint sets with $|V_1|=n_1$ and $|V_2|=n_2$, and let $E \subseteq V_1\times V_2$.  The graph $G = (V_1 \cup V_2,E)$ \index{graph!bipartite} \index{bipartite graph}is called a \emph{bipartite graph}.  Every edge in $G$ has its tail in $V_1$ and its head in $V_2$, and so 
\[
d_1 = \max\{ \outdeg_G(v_1) : v_1 \in V_1 \} \leq n_2
\]
and
\[
d_2 = \max\{ \indeg_G(v_2) : v_2 \in V_2 \} \leq n_1.
\]
Let $\alpha:V_1 \rightarrow \Gamma$ and $\beta: V_2 \rightarrow \Gamma$ be one-to-one functions from the vertices of $G$ to a group $\Gamma$.    We define
\[
\alpha(V_1) \overset{G}{+} \beta(V_2) = \{\alpha(v_1)+\beta(v_2) : (v_1,v_2) \in E\}.
\]
For all bipartite graphs $G$, we have 
\beq   \label{CH:ABineq}
|\alpha(V_1) \overset{G}{+} \beta(V_2)| \geq \max(d_1,d_2)
\eeq
for every group $\Gamma$ and all one-to-one maps $\alpha: V_1 \rightarrow \Gamma$ and $\beta: V_2 \rightarrow \Gamma$.

Consider the  \emph{complete bipartite graph}  $K_{n_1,n_2} = (V_1\cup V_2, V_1\times V_2)$. \index{complete bipartite graph} \index{graph!complete bipartite}  
We have $d_1 = \outdeg_G(v_1) = n_2$ for all $v_1 \in V_1$ and $d_2 = \indeg_G(v_2) = n_1$ for all $v_2 \in V_2$.
If $\Gamma = \Z$ and $\alpha: V_1 \rightarrow \Gamma$ and $\beta: V_2 \rightarrow \Gamma$ are one-to-one functions, 
then $|\alpha(V_1)  \overset{G}{+} \beta(V_2)|  \geq d_1+d_2-1$.   If $p$ is a prime and $\Gamma = \Z/p\Z$, then the Cauchy-Davenport theorem states that $|\alpha(V_1) \overset{G}{+} \beta(V_2)| \geq \min(d_1+d_2-1,p)$.   In particular, if  $\max(d_1,d_2)\leq p-1$ and $\min(d_1,d_2)\geq 2$, then $|\alpha(V_1) \overset{G}{+} \beta(V_2)| \geq \max(d_1,d_2)+1$.
One might guess that this is always a lower bound for $|\alpha(V_1) \overset{G}{+} \beta(V_2)|$, but the following beautiful construction by Josh Greene shows that inequality~\eqref{CH:ABineq} is best possible.

Let $A$ and $B$ be finite subsets of an abelian group $\Gamma$.  For every $x \in \Gamma$, we define the \emph{representation function} 
\[
r_{A,B}(x) = |\{(a,b) \in A\times B : a+b = x\}|.
\]
We construct the bipartite graph $G = (V_1\cup V_2,E)$, where
\[
V_1 = -A
\]
\[
V_2 = A+B
\]
and 
\[
E = \{ (-a,a+b) : b\in B\}.
\]
For all $v_1 \in V_1$ and $v_2 \in V_2$ we have
\[
\outdeg_G(v_1) = |B|  
\]
and 
\[
\indeg_G(v_2) = r_{A,B}(v_2).
\]
Then
\[
d_1 = \max\{\outdeg_G(v_1) : v_1 \in V_1\} = |B|
\]
and
\[
d_2 = \max\{\indeg_G(v_2) : v_2 \in V_2\} = \max\{ r_{A,B}(x) : x \in A+B \} \leq |B|.
\]
Define $\alpha:V_1 \rightarrow \Gamma$ by $\alpha(v_1) = v_1$ and 
 $\beta:V_2 \rightarrow \Gamma$ by $\beta(v_2) = v_2$.  Then
 \[
 \alpha(V_1) \overset{G}{+} \beta(V_2) = B
 \]
 and
 \[
 |\alpha(V_1) \overset{G}{+} \beta(V_2)| = |B| = \max(d_1,d_2).
 \]

Greene applied his construction in the following case.
Consider the Fermat prime $p = 257 = 2^{2^3} + 1$ and the finite field $\Gamma = \Z/257\Z$.  Let
\[
A = B = \{0\} \cup \{\pm 2^i : i = 0,1,\ldots,7\} \subseteq \Z/257\Z
\]
and 
\[
A+A = \{ a+a' : a,a' \in a\} \subseteq \Z/257\Z.
\]
Then $|A| = 17$ and $|A+A|=105$.   Note that every element of $A+A$ can be written as the sum of two distinct elements of $A$, since $0+0 = 1 + (-1)$, $2^7+2^7 = 0 + (-1)$, and $2^i + 2^i = 0 + 2^{i+1}$ for $i = 0,1,\ldots,6$.  Therefore,
$r_{A,A}(x) \geq 2$ for all $x \in A+A$, and so $d_1 = 105$ and $d_2 \geq 2$.

We conclude with a nice application of Sidon sets.  A \emph{Sidon set} is a subset $A$ of an abelian group such that every element of $A+A$ has a unique representation as the sum of two elements of $A$.  Equivalently,  $r_{A,A}(x) = 1$ for all $x \in A+A$.

Let $G = (V,E)$ be an undirected graph, and let $\alpha:V \rightarrow \Gamma$ and $\gamma: E \rightarrow \Gamma$ be one-to-one functions from the vertices and edges of $G$ into a group $\Gamma$.
Consider the set $\{\alpha(v)+ \gamma(e) : v \in e\}$.   If the maximum degree of a vertex in $V$ is $n$, then $|\{ \alpha(v) +\gamma(e) : v \in e\}| \geq n$.

\begin{theorem}[Jacob Fox]
Let $K_n$ denote the complete graph on $n$ vertices.  There are one-to-one functions $\alpha:V \rightarrow \Z$ and $\gamma:E \rightarrow \Z$ such that $|\{ \alpha(v)+ \gamma(e) : v \in e\}| = n$.
\end{theorem}

\begin{proof}
Denote the vertices of $K_n$ by $V = \{v_1,\ldots,v_n\}$ and the edges of $G$ by $E = \{ e_{i,j} = \{v_i,v_j\} : i,j = 1,\ldots,n\}$.  Let $A = \{a_1,\ldots,a_n\}$ be a Sidon set, and define the functions $\alpha$ and $\gamma$ by $\alpha(v_i) = -a_i$ for $i=1,\ldots,n$ and $\gamma(e_{i,j}) = a_i + a_j$ for $i,j =1,\ldots,n$.  Then $\{\alpha(v)+ \gamma(e) : v \in e\} = A$.  This completes the proof.
\end{proof}

\providecommand{\bysame}{\leavevmode\hbox to3em{\hrulefill}\thinspace}
\providecommand{\MR}{\relax\ifhmode\unskip\space\fi MR }
\providecommand{\MRhref}[2]{%
  \href{http://www.ams.org/mathscinet-getitem?mr=#1}{#2}
}
\providecommand{\href}[2]{#2}

\end{document}